\documentclass[12pt,a4paper]{article}
\usepackage{t1enc}
\usepackage[latin1]{inputenc}
\usepackage[english]{babel}
\usepackage{amssymb}
\usepackage{amsmath}
\usepackage{graphics}

\usepackage{amsthm, amssymb}
\usepackage{amsfonts}
\usepackage{epsfig,multicol}

\newtheorem{theorem}[subsubsection]{Theorem}
\newtheorem{proposition}[subsubsection]{Proposition}
\newtheorem{lemma}[subsubsection]{Lemma}
\newtheorem{corollary}[subsubsection]{Corollary}
\newtheorem{definition}[subsubsection]{Definition}

\newtheorem{remark}[subsubsection]{Remark}

\newfont{\gothic} { ygoth scaled \magstep{1.5}}

\newcommand{\5}{\vskip 5pt}

\def\<{\langle}
\def\>{\rangle}

\begin{document}

\def\hpic #1 #2 {\mbox{$\begin{array}[c]{l} \epsfig{file=#1,height=#2}
\end{array}$}}
 
\def\vpic #1 #2 {\mbox{$\begin{array}[c]{l} \epsfig{file=#1,width=#2}
\end{array}$}}

\title{Forked Temperley-Lieb Algebras and Intermediate Subfactors}


\author{Pinhas Grossman}




\maketitle


\def\vpic #1 #2 {\mbox{$\begin{array}[c]{l} \includegraphics[width=#2]{#1}
\end{array}$}}
\begin{abstract}
We consider noncommuting pairs $P,Q$ of intermediate subfactors of an irreducible, finite-index inclusion $N \subset M$ of II$_1$ factors such that $P$ and $Q$ are supertransitive with Jones index less than $4$ over $N$. We show that up to isomorphism of the standard invariant, there is a unique such pair corresponding to each even value $[P:N]=4cos^2\frac{\pi}{2n}$ but none for the odd values $[P:N]=4cos^2 \frac{\pi}{2n+1} $. 

We also classify the angle values which occur between pairs of intermediate subfactors with small index over their intersection: if \\$[P:N],[Q:N] < 4$, then the unique nontrivial angle value is always $cos^{-1}\displaystyle \frac{1}{[P:N]-1} $.
\end{abstract}
\section{Introduction}
A fundamental example of a subfactor is the fixed-point algebra of an outer action of a finite group on a von Neumann algebra with trivial center. In this case the structure of the subfactor is determined by the structure of the group. One then thinks of a general subfactor as a ``quantum'' version of a finite group, and subfactor theory as a ``non-commutative Galois theory''. In this spirit it is natural to consider an intermediate subfactor $N \subset P \subset M$ as an analogue of a subgroup, and indeed, the intermediate subfactors of the fixed-point subfactor of an action of a finite group are precisely the fixed-point subfactors of its subgroups.

The main classification theory for subfactors involves a \textit{standard invariant}, which was axiomatized by Popa \cite{P20} and then described by Jones using a diagrammatic apparatus called \textit{planar algebras} \cite{J18}. A major structural feature of subfactors is duality: there are two algebraic structures applicable to the standard invariant (analagous to Hopf-algebra duality for groups) which are reflected in the planar algebra picture by a choice of shading. Bisch discovered the remarkable fact that intermediate subfactors are characterized by ``biprojections''- elements of the standard invariant that are projections in both algebra structures \cite{Bs3}. Bisch's theorem thus reduces the problem of intermediate subfactors to a problem of planar algebras.

Bisch and Jones then studied the planar algebra generated by a single biprojection, which yields a generic construction of an intermediate subfactor, or more generally, a chain of intermediate subfactors \cite{BJ}. This led naturally to the approach of trying to construct more general families of intermediate subfactors by considering planar algebras generated by multiple biprojections.  

There is a notion of commutativity for pairs of intermediate subfactors (defined simply by commutativity of the corresponding projections in $L^2(M)$.) Sano and Watatani considered the \textit{set of angles} between two subfactors, a numerical invariant which measures the degree of noncommutativity of a pair of intermediate subfactors \cite{SaW}. It turns out that pairs of commuting subfactors can be constructed simply via a tensor product, but the construction of intermediate subfactors with nontrivial angles has proven much more difficult.

In \cite{GJ}, Jones and the present author set out to construct generic pairs of noncommuting intermediate subfactors by assuming \textit{no extra structure}, i.e. that the standard invariants of the elementary subfactors involved were just Temperley-Lieb algebras, a situation also referred to as \textit{supertransitivity}. It was hoped that since the bimodule components of the intermediate subfactors obey the Temperley-Lieb fusion rules, that information together with invariants such as angles and indices would provide sufficient rigidity to suggest appropriate planar relations among biprojections and yield generic constructions of planar algebras generated by two biprojections. 

Instead we found the surprising result that there are essentially only two quadrilaterals with no extra structure. One is the fixed point algebra of an outer action of  $S_3$ on a factor (where the intermediate subfactors are fixed point algebras of distinct order $2$ subgroups) and the other is a new example with irrational angles and indices coming from the GHJ family of subfactors. 

The dearth of examples with no extra structure does not however mean that the original project of constructing planar algebras generated by multiple biprojections is unsound. Rather, because of the rigidity of intermediate subfactors, we must adjust our hypotheses if we wish to obtain more examples.

In the present work we continue this approach by considering a (noncommuting, irreducible, finite index) quadrilateral \vpic{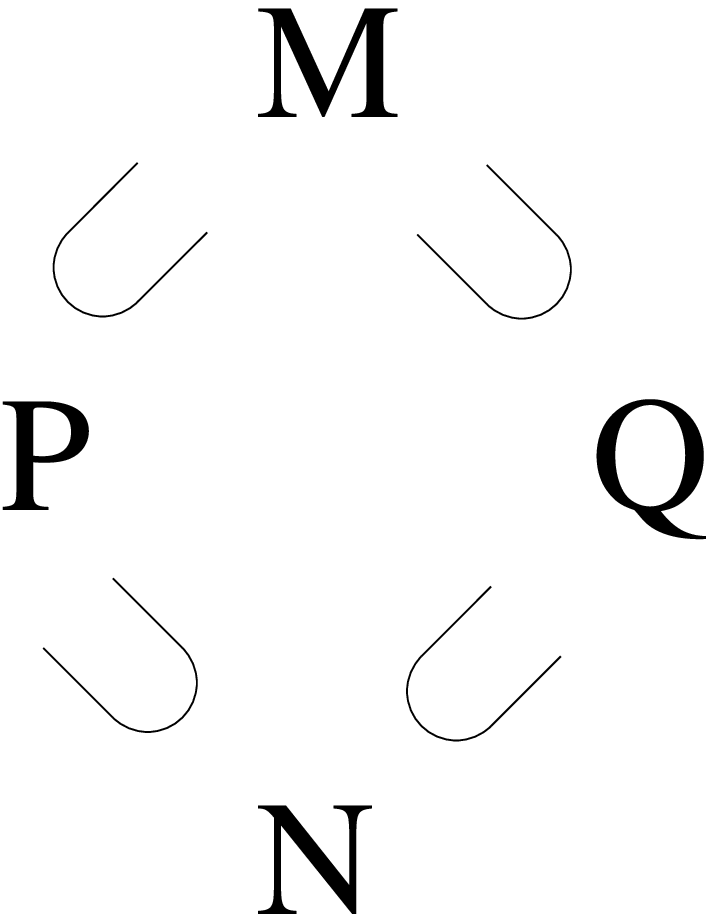} {0.6 in}  in which only the lower two elementary subfactors $N \subset P$ and $N \subset Q $ are assumed to be supertransitive. We also assume that the lower indices $[P:N] $ and $[Q:N]$ are less than $4$ (and therefore equal to $4cos^2\frac{\pi}{k} $ for some $k \geq 3$ by Jones' celebrated theorem in \cite{J3}).

We first show that for any quadrilateral with lower indices $[P:N] $ and $[Q:N]$ less than $4$ there is a unique nontrivial angle value, equal to $cos^{-1}\displaystyle \frac{1}{[P:N]-1} $. We then describe a series of examples in the GHJ family-first discovered in \cite{GJ}- of quadrilaterals whose lower subfactors are supertransitive with indices $[P:N]=[Q:N]=4cos^2\frac{\pi}{2k}, k\geq 3 $. 

Ultimately we prove the following result:

\begin{theorem}
The series of GHJ subfactors for $D_n$ at the trivalent vertex gives a complete list of noncommuting irreducible hyperfinite quadrilaterals such that $N \subset P$ and $N \subset Q$ are hyperfinite supertransitive subfactors with index less than $4$.

\end{theorem}

In particular it turns out there is a unique such quadrilateral for each even value $[P:N]=4cos^2\frac{\pi}{2k} $ but none for the odd values $[P:N]=4cos^2\frac{\pi}{2k-1} $.

The methods of proof involve both the planar apparatus for intermediate subfactors developed in \cite{GJ} and a notion of \textit{forked Temperley-Lieb algebras}, in which a sequence of Temperley-Lieb projections $e_1,e_2,e_3,...$ in a von Neumann algebra can be initially extended by either of two orthogonal projections $p$ and $q$. This is a special case of $AF$ algebras associated to $T$-shaped graphs, which were studied by Evans and Gould in \cite{EG}. It turns out this situation precisely captures the structure of noncommuting supertransitive intermediate subfactors with small index. 

\textbf{Acknowledgement.} I would like to thank my advisor, Vaughan Jones, for his support and encouragement, and for many useful conversations regarding this work.

\section{Background}
\subsection{Subfactors and planar algebras}
The best reference for the elementary theory of subfactors is Jones' original paper (\cite{J3}). We review the basics here.

A \textit{II$_{1}$ factor} $M$ is an infinite dimensional von Neumann algebra with trivial center which admits a faithful positive-definite trace. Such a trace, normalized to equal one at the identity, is necessarily unique. The trace induces an inner product on $M$, and its Hilbert space completion is denoted by $L^2(M)$. $M$ acts on $L^2(M)$ by left (or right) multiplication, continuously extended.

 Let $N \subset M$ be a unital inclusion of II$_1$ factors. Then $L^{2}(N)$ can be
 identified with a subspace of $L^{2}(M)$. Let $e_{1}$ denote the corresponding projection on $L^{2}(M)$, and let $M_{1}$ be the von Neumann algebra generated by $M$ (thought of as an algebra of left multiplication operators) and $e_{1}$. This procedure is called the \textit{basic construction}. 

If $M$ is finitely generated as an $N$-module, then $M_1$ is a II$_1$ factor. In this case the index $[M:N]$ is defined to equal $\tau^{-1} $, where $\tau=tr(e_1)$. Otherwise the index is defined to be infinity. Some properties of the index are: $[M:N] \geq 1 $, $[M_1:M]=[M:N] $, and if $N \subseteq P \subseteq M$, then $[M:N]=[M:P][P:N] $.   

It turns out that not all numbers greater than $1$ can occur as index values, by the following remarkable theorem of Jones.

\begin{theorem}\textbf{(Jones)}
Let $N \subseteq M$ be a unital inclusion of II$_1$ factors with finite index. Then either $[M:N] \geq 4$ or there is an integer $k \geq 3$ such that $[M:N]=4cos^2 \frac{\pi}{k}$.
\end{theorem}

If $[M:N]<\infty $, one can iterate the basic construction to obtain a sequence of projections $e_{1},e_{2}...$, and 
a tower of II$_1$ factors $M_{-1} \subset M_{0} \subset M_{1} \subset M_{2} \subset ...$,
 where $M_{-1}=N$, $M_{0}=M$, $e_{k}$ is the projection onto
 $L^{2}(M_{k-2})$ in $\textbf{B}(L^{2}(M_{k-1}))$, and $M_{k}$ is the von Neumann algebra generated
 by $M_{k-1}$ and $e_{k}$, for $k \geq 1$. Restricting the tower to those elements which 
commute with $N$, one obtains a smaller tower of algebras- which turn out to be finite-dimensional- called the tower of relative
 commutants $N'\cap M_k$. 

The tower of relative commutants can also be thought of as algebras of bimodule intertwiners of tensor powers of the $N-N$ bimodule $L^2(M)$.  

\begin{proposition}
As $N-N$ bimodules, $L^2(M_{k}) \cong L^2(M) \otimes_{N} ... \otimes_{N} L^2(M)$, ($k+1$ factors). Moreover, $Hom_{N-N} L^2(M) \cong N' \cap M_{2k+1}$. So an $N-N$ bimodule decomposition of $\otimes_N^{k+1} L^2(M)$ corresponds to
 a decomposition of the identity in $N' \cap M_{2k+1}$. Under this correspondence
projections in $N' \cap M_{2k+1}$ correspond to submodules of $\otimes_N^{k+1}  L^2(M)$, minimal projections 
correspond to irreducible submodules (those which have no proper nonzero closed submodules),
 and minimal central projections to equivalence classes of irreducible submodules.
\end{proposition}

\begin{remark}
The even relative commutants correspond to intertwiners of ``mixed'' bimodules over $N$ and $M$.
\end{remark}

The lattice of finite-dimensional von Neumann algebras\\\\
$\begin{array}{ccccccccc}
N' \cap N &\subseteq &N' \cap M &\subseteq &N' \cap M_1 &\subseteq &N' \cap M_2 &\subseteq &... \cr
 & & \cup & & \cup &  & \cup & & \cr

  &  &M' \cap M & \subseteq & M' \cap M_1 & \subseteq & M' \cap M_2 & \subseteq & ...
\end{array}$\\\\

is called the \textit{standard invariant} and was axiomatized by Popa \cite{P20}. Popa showed that the standard invariant is a complete invariant under certain conditions called \textit{strong amenability}, and also introduced a construction of subfactors having a given standard invariant.

Note that the standard invariant always contains the Jones projections $e_1,e_2,e_3,...$, but in general will contain more structure as well. The Jones projections satisfy the following relations: (i) $e_i e_{i\pm1} e_i=\tau e_i$; and
(ii) $e_i e_j=e_j e_i$ if $|i-j| \geq 2$, where $\tau=[M:N]^{-1} $. The algebras they generate are in a sense ``minimal'' standard invariants, and are known as the Temperley-Lieb algebras after their original appearance in \cite{TL}. 

In \cite{Kff2} Kauffman introduced a pictorial representation of the Temperley-Lieb algebras related to the skein theory of knot diagrams. It was soon noticed by Jones and others that many seemingly complicated linear relations in the standard invariants of subfactors had simple pictorial representations. To capture this remarkable phenomenon in a rigorous and systematic fashion, Jones introduced the notion of a \textit{planar algebra}.

A \textit{planar tangle} consists of a disc in the plane, together with some disjoint internal discs and nonintersecting strings connecting certain boundary points of the discs, all defined up to planar isotopy. There is a natural notion of composition of tangles (by isotoping one tangle to fit in an internal disc of another tangle and removing the boundary). One then thinks of the internal discs as ``inputs'' and the external boundary as the ``output'', and the collection of planar tangles forms an operad.

If $N \subset M$ is a subfactor, its planar algebra is essentially an action of the operad of planar tangles on the central vectors of the tensor powers of the $N-N$ bimodule $M$, i.e. the standard invariant. In this picture, the trivial tangles-those with no internal discs- correspond to the Temperley-Lieb algebras, which are present for any subfactor, and more complicated tangles may reflect additional structure. In this way, the algebraic and combinatorial structure of a subfactor is captured in terms of the geometry of the plane.

Jones' axioms for a planar algebra turned out to be exactly analogous to Popa's axioms for a $\lambda$-lattice in his description of the standard invariant, so by Popa's theorem a subfactor can actually be constructed by describing its planar algebra, as in the work of Bisch and Jones (\cite{BJ}).

For more precise details on definitions and examples of planar algebras and the meaning of various pictures we direct the reader to \cite{J18}. Note, however, that we adopt the drawing conventions of \cite{GJ}.  

\subsection{Intermediate subfactors}
We now recall some facts about intermediate subfactors; for more details and proofs we refer the reader to \cite{GJ}.

Let $N \subset M$ be an irreducible inclusion of $II_{1}$ factors with finite index.  (Irreducible here means that $N' \cap M = \mathbb{C}Id$, or equivalently that $M$ is an irreducible $N-M $ bimodule). Consider also the dual  inclusion from the basic construction $M \subset M_{1}$. The first relative commutants $N' \cap M_{1}$ and $M' \cap M_{2}$ have the same vector space dimension, and the map $\phi:$ \vpic{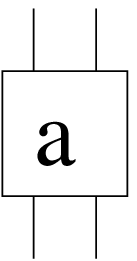} {0.3 in}  $\mapsto$ \vpic{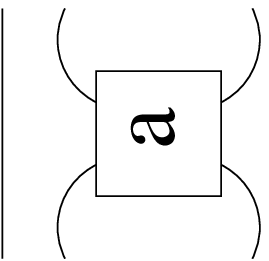} {0.6 in}  is a linear isomorphism.   
 Pulling back the multiplication from $M' \cap M_2$ via $\phi $ gives a second multiplication structure on $N' \cap M_1$, called ``comultiplication'' and denoted by the symbol $\circ$: if $a$ and $b$ are elements of $N' \cap M_{1}$, then 
$a \circ b =$ \vpic{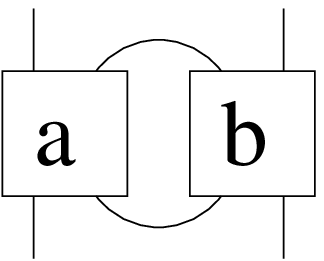} {0.7 in}  . The following result of Bisch characterizes intermediate subfactors in terms of the standard invariant.

\begin{theorem}\textbf{(Bisch)}
Let $p=p^*=p^2 $ be a projection in $N' \cap M_1 $. Then $p $ projects onto an intermediate subfactor $N \subseteq P \subseteq M $ iff $p$ is (a scalar multiple of) a projection in the dual algebra structure. 
\end{theorem}

Therefore an intermediate subfactor projection is called a \textit{biprojection}. Landau discovered the following important relation among biprojections.
 
\begin{theorem}\label{landauproj}
\textbf{(Landau)} Let $P$ and $Q$ be intermediate subfactors with biprojections $e_P$ and $e_Q$. Then $e_{P} \circ e_{Q}=$  \vpic{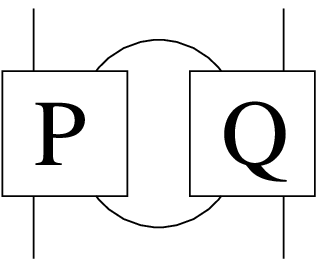}  {0.6 in}  $=\delta tr (e_{P}e_{Q}) e_{PQ}$, where $e_{PQ}$ is the projection onto the vector space $PQ \subseteq M$.
\end{theorem}
The following trace formulas follow. 
\begin{proposition} \label{mulfor}
We have $tr(e_{PQ}) tr(e_{P}e_{Q}) = tr(e_{P})tr(e_{Q})$. Also,
$tr(e_{P}e_{Q}) = \displaystyle \frac{1}{dim_{M}L^{2}(\bar{P}\bar{Q})}$, where $\bar{P} $ and $\bar{Q} $ are the dual factors from the basic construction for $P \subset M $ and $Q \subset M$ respectively.
\end{proposition}

A pair of intermediate subfactors of $M$ is said to commute if the corresponding biprojections commute in the planar algebra. The subfactors are said to cocommute if the biprojections cocommute. A quadrilateral of II$_1$ factors is a diagram \vpic{quad} {0.6 in}  such that $P \vee Q =M $ and $P \wedge Q=N$. The following result provides a useful characterization of commutativity and cocommutativity.

\begin{theorem}
Let \vpic{quad} {0.6 in}  be a quadrilateral of $II_{1}$ factors,
 where $N \subset M$ is an irreducible finite-index inclusion. Consider 
the multiplication map from the bimodule tensor product $P \otimes_{N} Q$ to $M$. 
The quadrilateral commutes iff this map is injective and cocommutes iff the map is surjective.
\end{theorem}

\subsection{Supertransitive Subfactors}
We recall some facts about supertrnsitive subfactors, which were defined by Jones in \cite{J3}. Let $N \subset M$ be an inclusion of II$_{1}$ factors with associated tower 
$M_{-1} \subset M_{0} \subset M_{1} \subset ...$, where $M_{-1}=N$, $M_{0}=M$,
 and $M_{k+1}$, $k \geq 0$ is the von Neumann algebra on $L^{2}(M_{k})$ generated by
 $M_{k}$ and $e_{k+1}$, the projection onto $L^{2}(M_{k-1})$.
 Each $e_{k}$ commutes with $N$, so $\{ 1,e_{1},..,e_{k} \}$ generates a *-subalgebra, 
which we will call $TL_{k+1}$, of the $k^{th}$ relative commutant $N' \cap M_{k}$.

A finite-index subfactor $N \subseteq M$ is called $k$-supertransitive 
 $k$-{\rm supertransitive} (for $k>1$) if $N' \cap M_{k-1} = TL_k$. 
We will say $N\subseteq M$ is {\rm supertransitive} if it is $k$-supertransitive
for all $k$. Note that $N\subseteq M$ is $1$-supertransitive iff it is irreducible, i.e. $N'\cap M \cong \mathbb C$ and is $2$-supertransitive iff the $N-N$ bimodule $L^2(M)$ has two irreducible components. Supertransitivity of $N\subseteq M$ is the same as saying its principal graph 
is $A_n$ for some $n=2,3,4,..., \infty$.

\begin{lemma} \label{fusionrules}
Suppose  $N \subset M$ is supertransitive. If $[M:N] \geq 4$ then there is a sequence 
of irreducible $N-N$ bimodules
 $ V_{0}, V_{1}, V_{2}... $ such that $L^{2}(N) \cong V_{0}$,
 $L^{2}(M) \cong V_{0} \oplus V_{1}$, and $V_{i} \otimes V_{j} \cong \oplus_{k=|i-j|}^{i+j}V_{k}$.
 If $[M:N]=4cos^{2}(\frac{\pi}{n})$ then the sequence terminates at $V_{l}$, 
where $l=[\frac{n-2}{2}]$, and the fusion rule
 is: $V_{i} \otimes V_{j} \cong \oplus_{k=|i-j|}^{(\frac{n-2}{2})-|(\frac{n-2}{2})-(i+j)|}V_{k}$. (see \cite{BJ2} ) \\
In either case, we have $\displaystyle \dim_{N}V_{k}=[M:N]^{k}T_{2k+1}(\frac{1}{[M:N]})$, 
where $\{ T_{k}(x) \}$ is the sequence of polynomials defined recursively by
 $T_{0}(x)=0$, $p_{1}(x)=1$, and $T_{k+2}(x)=T_{k+1}(x)-xT_{k}(x)$. In particular, 
$dim_{N}V_{1}=[M:N]-1$ and $dim_{N}V_{2}=[M:N]^{2}-3[M:N]+1$.
\end{lemma}

\begin{remark} \label{modfusionrules}
If $N \subset M$ is $2k$-supertransitive, then there is a sequence 
of irreducible bimodules $V_{0}, ..., V_{k}$ for which the above fusion rules and 
dimension formula hold as long as $i + j \leq k$. 
\end{remark}

Now let \vpic{quad} {0.6 in}  be a quadrilateral of finite index subfactors. The four subfactors $N\subseteq P$,$N \subseteq Q$,$P \subseteq M$, and $Q\subseteq M$ are called the {\it elementary subfactors}. A quadrilateral is said to have {\rm no extra structure} if all the elementary subfactors are supertransitive. The following result from \cite{GJ} classifies noncommuting irreducible quadrilaterals with no extra structure.

\begin{theorem}\label{maintheorem}
Suppose \vpic{quad} {0.6 in}  
 is a quadrilateral of
subfactors with $N'\cap M =\mathbb C$, $[M:N]<\infty$ and no extra structure. 
Then either the quadrilateral commutes or  one
of the following two cases occurs:

a) $[M:N]=6$ and $N$ is the fixed point algebra for an outer action of $S_3$
on $M$ with $P$ and $Q$ being the fixed point algebras for two transpositions
in $S_3$. In this case the angle between $P$ and $Q$ is $\pi/3$.

b) $[M:N]=(2+\sqrt 2)^2$ and the planar algebra of $N\subseteq M$ is the same as
that coming from the GHJ subfactor (see \cite{GHJ}) constructed from the Coxeter graph $D_5$ with the distinguished vertex being the trivalent one. 
Each of the intermediate inclusions has index $2+\sqrt 2$ and the angle between
$P$ and $Q$ is $\theta=\cos^{-1} (\sqrt 2 - 1)$.
\end{theorem} 

\section{Supertransitive intermediate subfactors in the GHJ family}
\subsection{GHJ subfactors}

We recall the construction of the GHJ subfactors from the book of Goodman, de la Harpe and Jones \cite{GHJ}. Let $G$ be a Coxeter-Dynkin diagram of type $A$, $D$, or $E$, with a distinguished vertex $*$ and a bipartite structure. Let $A_0 \subset A_1$ be an inclusion of finite-dimensional von Neumann algebras the underlying graph of whose Bratteli diagram is $G$, with $*$ corresponding to a particular simple direct summand of $A_0$. 

Perform the basic construction on the inclusion $A_0 \subset A_1$ with respect to the Markov trace, and iterate to obtain a tower $A_0 \subset A_1 \subset A_2 \subset A_3 \subset ...$ . Inside this tower are the Jones projections $e_1, e_2, e_3, ...$. Let $B_i$ be the subalgebra of $A_i$ generated by $e_1,...,e_{i-1}$, for each $i=2,3,...$. Then the towers $B_i \subset A_i$ form commuting squares, and the von Neumann algebra $B$ generated by $\cup_{i=2}^{\infty}B_i $ is a subfactor of the von Neumann algebra $A$ generated by $\cup_{i=2}^{\infty}A_i $ on $L^2(\cup_{i=2}^{\infty}A_i ) $. This subfactor has finite index but is not in general irreducible. Let $r$ be the projection in $A_0$ corresponding to $*$. Then $r \in B'$, and $rB \subset rAr $ is an irreducible subfactor with finite index, called the GHJ subfactor for $G$, $*$. 

\begin{remark} 
The Jones projections $e_1,e_2,e_3,...$ satisfy the Temperley-Lieb relations:\\
(i) $e_i e_{i\pm1} e_i=\tau e_i$\\
(ii) $e_i e_j=e_j e_i$ if $|i-j| \geq 2$\\
(iii) $tr(e_i w)=\tau tr(w) $ if $w$ is a word on $e_1,...,e_{i-1}$.\\
Here $\tau=\displaystyle \frac{1}{4cos^2 \frac{\pi}{k}}$, where $k$ is an integer associated to $G$ called \textit{the Coxeter number}. For each $n$, the Coxeter number of $A_n$ is $n+1$, the Coxeter number of $D_n$ is $2n-2$, and the Coxeter numbers of $E_6$, $E_7$, and $E_8$ are $12$, $18$, and $30$ respectively.

The cutdown Jones projections $re_1,re_2,re_3,...$ satisfy the same relations, and we will suppress the ``$r$'' when dealing with them.  
\end{remark}

\begin{remark} \label{braid_rep}
Because of the Temperley-Lieb relations, there is a unitary representation of the braid group inside the algebras $B_i$ (or $rB_i$), sending the usual braid group generators $\sigma_i$ to $g_i=(t+1)e_i-1$, where $t=e^{2\pi i/k}$, $k$ again being the Coxeter number of $G$.
\end{remark}

The principal graphs of the GHJ subfactors were computed by Okamoto in \cite{Ok}.

\subsection{The GHJ subfactor for $D_n$ at the trivalent vertex, and a pair of intermediate subfactors} \label{ghj_dn}

The Bratteli diagram for the tower $rA_0r \subset rA_1r \subset rA_2r \subset ...$ is determined by the \textit{string algebra} for $G$ based at $*$ (see \cite{EK7}). Thus for example the initial steps $rA_0r \subset rA_1r \subset rA_2r \subset rA_3r$ for $D_n$ at the trivalent vertex have the following diagram:

 \vpic{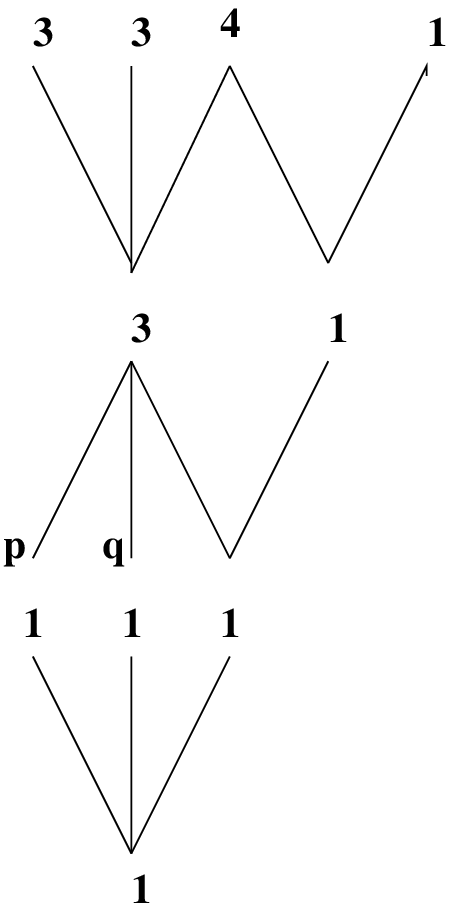} {1.35 in}  

Let $N \subset M$ be the GHJ subfactor for $D_n$ at the trivalent vertex. Let $p$ and $q$ be the projections corresponding to the two terminal vertices of the fork, as in the above figure. Then evidently $p$ and $q$ are orthogonal, and commute with $e_2,e_3,...$. From the string algebra construction it follows that $p$ and $q$ have the same trace as the $e_i$, and therefore $pe_1p=\tau p$ and $e_1 p e_1=\tau e_1$, and similarly for $q$. Thus $p$ and $q$ can each be used independently to initially extend the sequence of Temperley-Lieb projections which generate $N$.
We recall the following theorem from \cite{J3}:
\begin{theorem}\textbf{(Jones)}
Let $M$ be a von Neumann algebra generated by a sequence of projections $e_1,e_2,e_3,...$ satisfying the Temperley-Lieb relations with $\tau=\displaystyle \frac{1}{4cos^2\frac{\pi}{k}}$, and let $N$ be the subalgebra generated by $e_2,e_3,...$. Then $N$ and $M$ are II$_1$ factors and the principal graph for $N \subset M$ is $A_{k-1}$.
\end{theorem}

Let $P$ be the von Neumann algebra generated by $N$ and $p$, and let $Q$ be the von Neumann algebra generated by $N$ and $q$. Because $p$ extends initially the sequence of Temperley-Lieb projections which generate $N$, the above theorem applies with $\tau=\displaystyle \frac{1}{4cos^2 \frac{\pi}{2n-2}}$, since the Coxeter number for $D_n$ is $2n-2$. Similarly for $N \subset Q$. Therefore the principal graphs for the subfactors $N \subset P$ and $N \subset Q$ are both $A_{2n-3}$.

\begin{remark}\label{2reps}
In a similar way, the braid group representation of \ref{braid_rep} can be initially extended in two different ways, depending on whether one chooses $p$ or $q$ as $e_0$.
\end{remark}

\subsection{Example: The $D_5$ case, and intermediate subfactors with no extra structure}

For $D_5$, the GHJ subfactor at the trivalent vertex has index $6 + 4 \sqrt{2}=(2+\sqrt{2})^2$. Moreover, the Coxeter number for $D_5$ is $8$, and  $4cos^2(\frac{\pi}{8})=2+\sqrt{2} $. So we have $[M:P]=[M:Q]=[P:N]=[Q:N]=2+\sqrt{2}$. This was one construction of the quadrilateral with no extra structure in \cite{GJ}. The original construction involved the ``GHJ pair'' $\tilde{P}$ and $\tilde{Q}$ coming from the two braid group representations of \ref{2reps}. The GHJ pair also has index $2+\sqrt{2}$ in $M$. One might have expected these two pairs of intermediate subfactors to coincide in this case (in general though they cannot coincide since although we always have $[M:\tilde{P}]=[P:N] $, the other indices are usually different.)  

But in fact the two pairs are distinct, and we get two distinct quadrilaterals of total index $(2+\sqrt{2})^2 $. Neither quadrilateral commutes nor cocommutes, but they do commute with each other. And they are also isomorphic as quadrilaterals. There are also intermediate subfactors $R$ and $S$ of index and co-index $2$, coming from a period two automorphism exchanging $P$ and $Q$ (and $\tilde{P}$ and $\tilde{Q}$). It can be shown that there are no other intermediate subfactors, and so the GHJ subfactor for $D_5$ at the trivalent vertex has the intermediate subfactor lattice computed in \cite{GJ}:

\vpic{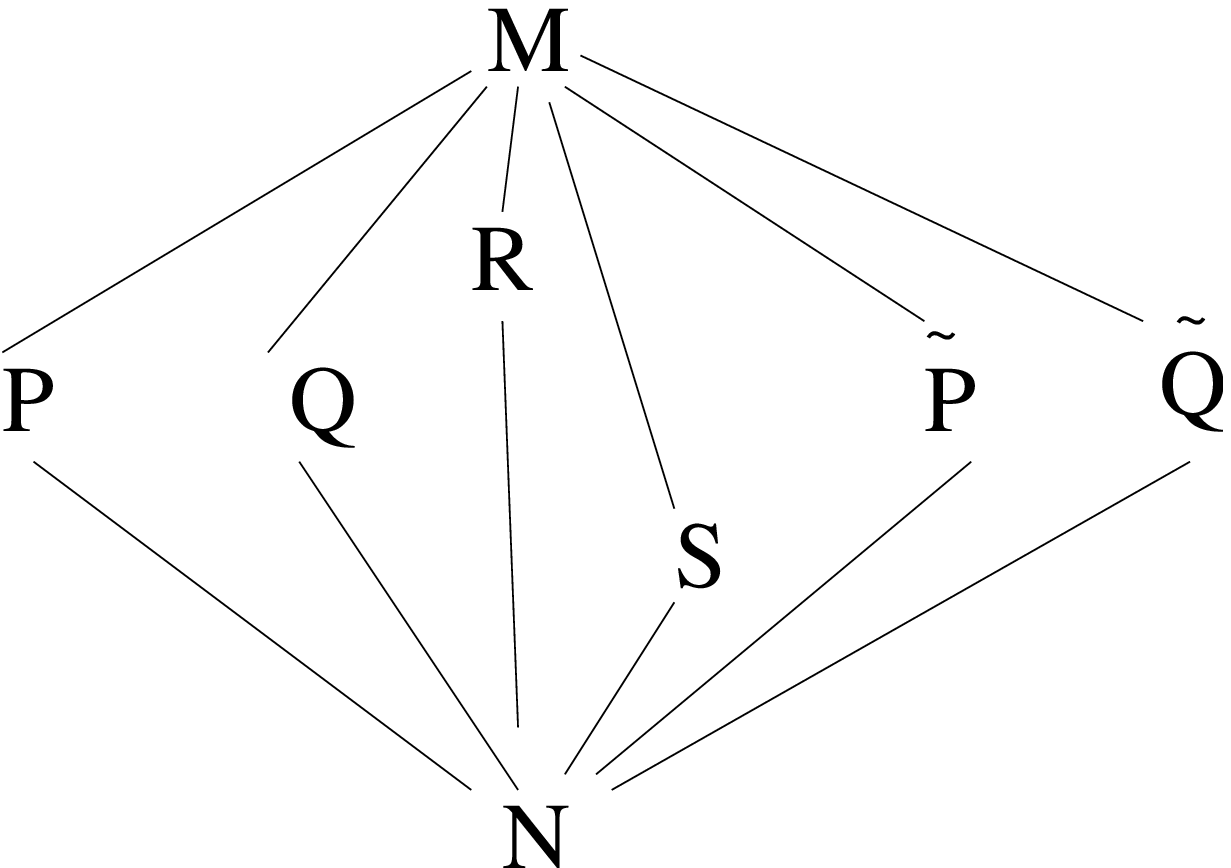} {2.0 in}

\section{Angles between intermediate subfactors with small index}
\subsection{Determination of possible angle values}
In \cite{SaW}, Sano and Watatani studied the notion of the set of angles between subfactors. We recall the definition.

\begin{definition}
Let $P,Q$ be subfactors of $M$. The \textit{set of angles} between $P$ and $Q$ is given by the spectrum of $cos^{-1} \sqrt{e_P e_Q e_P} $. 
\end{definition}

It turns out that for pairs of intermediate subfactors with index less than $4$ over the intersection, the angles are completely determined by the index.

We first recall the following result from \cite{PP2}.

\begin{lemma} \label{pplemma} (\textbf{Pimsner-Popa})
If the $N-N$-bimodule decomposition of $L^{2}(M)$ contains $k$ copies of the
 $N-N$-bimodule $R$, then $k \leq dim_{N} R$. In particular, $L^{2}(M)$ contains only 
one copy of $L^{2}(N)$.
\end{lemma}

We will also need some results from \cite{GJ}.

\begin{lemma} \label{pisoq}
If $N\subseteq P$ and $N\subseteq Q$ are 2-transitive  and the quadrilateral does not 
commute then $L^{2}(P) \cong L^{2}(Q)$ as $N-N$-bimodules, and therefore $[P:N]=[Q:N]$.
\end{lemma}

\begin{lemma}\label{lambda_calc}
If $N \subset P$ and $N \subset Q$ are $2$-transitive, then 
$e_{P}e_{Q}e_{P}=e_{N} + \lambda (e_{P}-e_{N})$, 
where $\lambda=\displaystyle \frac{tr(e_{\bar{P}\bar{Q}})^{-1}-1}{[P:N]-1}$. 
\end{lemma}

\begin{lemma} \label{pqdecomp}
If $N\subseteq P$ and $N\subseteq Q$ are 4-supertransitive and the quadrilateral does not commute
 then the $N-N$-bimodule $L^{2}(PQ)$ isomorphic to one of the 
following: $V_{0} \oplus 2V_{1} \oplus V_{2}$, $V_{0} \oplus 3V_{1} \oplus V_{2}$, 
or $V_{0} \oplus 3V_{1}$, where the $V_{i}$ are as in \ref{fusionrules}   
(for the 4-supertransitive inclusion $N \subset P$).
\end{lemma}
\begin{theorem}\label{angle_thm} 
Let \vpic{quad} {0.6 in}  be a noncommuting irreducible quadrilateral with finite index, and suppose that $[P:N]$ and $[Q:N]$ are less than $4$. Then $Ang(P,Q)$ contains a unique nontrivial value, equal to $\cos^{-1} \displaystyle \frac{1}{[P:N]-1}$.
\end{theorem}

\begin{proof}
Since $[P:N]$ and $[Q:N]$ are less than $4$, the inclusions $N \subset P$ and $N \subset Q$ have principal graphs equal to $A_n$ or $D_{2n}$ for some $n$, or $E_6$ or $E_8$. The only candidate which lacks $2$-supertransitivity is $D_4$- however this is impossible by the noncommutativity hypothesis, since the $N-N$ bimodules of a $D_4$ subfactor have dimension $1$, and the corresponding projections would necessarily be central by \ref{pplemma}. So we may assume $2$-supertransitivity of $N \subset P$ and $N \subset Q$. Then by \ref{pisoq} $L^2(P) \cong L^2(Q)$ as $N-N$ bimodules.

By $2$-supertransitivity and \ref{lambda_calc} we have $e_{P}e_{Q}e_{P}=e_{N} + \lambda (e_{P}-e_{N})$, 
where $\lambda=\displaystyle \frac{tr(e_{\bar{P}\bar{Q}})^{-1}-1}{[P:N]-1}$, so there is a unique nontrivial angle value, and it is equal to $cos^{-1} \sqrt{\frac{tr(e_{\bar{P}\bar{Q}})^{-1}-1}{[P:N]-1}} $, which by \ref{mulfor} can be rewritten as \\ $cos^{-1} \sqrt{\displaystyle \frac{\frac{[M:N]}{dim_N L^2(\bar{P}\bar{Q})}-1}{[P:N]-1}} =cos^{-1} \sqrt{\displaystyle \frac{tr(e_P e_Q)[M:N]-1}{[P:N]-1}}= \\ cos^{-1} \sqrt{\displaystyle \frac{\frac{tr(e_P)tr(e_Q)[M:N]}{tr(e_{PQ})}-1}{[P:N]-1}}=
cos^{-1} \sqrt{\displaystyle \frac{\frac{[M:P]^{-2}[M:N]^{2}}{[M:N]tr(e_{PQ})}-1}{[P:N]-1}} = \\ cos^{-1} \sqrt{\displaystyle \frac{\frac{[P:N]^2}{dim_N L^2(PQ)}-1}{[P:N]-1}} $.

Suppose that $N \subset P$ is also $4$-supertransitive. Then by \ref{pqdecomp} we have that the $N-N$-bimodule $L^{2}(PQ)$ isomorphic to one of the following: $V_{0} \oplus 2V_{1} \oplus V_{2}$, $V_{0} \oplus 3V_{1} \oplus V_{2}$, or $V_{0} \oplus 3V_{1}$, where $dim_N V_0=1 $, $dim_N V_1=[P:N]-1 $, and $dim_N V_2=$\\  $[P:N]^2-3[P:N]+1$. However, since $[P:N]<4 $, by \ref{pplemma} $L^2(M)$ can contain at most $2$ copies of $V_1$. So in fact $L^2(PQ) \cong V_0 \oplus 2V_1 \oplus V_2$, and $dim_N L^2(PQ)=1+2([P:N]-1)+[P:N]^2-3[P:N]+1=[P:N]^2-[P:N]=[P:N]([P:N]-1)$.

Plugging this into the angle formula above, we find that the angle is \\$cos^{-1} \sqrt{\displaystyle \frac{\frac{[P:N]^2}{[P:N]([P:N]-1)}-1}{[P:N]-1}} =cos^{-1} \sqrt{\displaystyle \frac{ \frac{1}{[P:N]-1} } {[P:N]-1} } =cos^{-1}\displaystyle \frac{1}{[P:N]-1}$.

All this assumed $4$-supertransitivity. However, the only admissible principal graphs lacking $4$-supertransitivity are $D_6$ and $E_6$. But a noncommuting quadrilateral whose lower subfactors have principal graph $E_6$ cannot exist (see \cite{G}). For any noncommuting quadrilateral whose lower subfactors have principal graph $D_6$, by the $D_6$ fusion rules we have $L^2(PQ)=L^2(M) \cong L^2(N) \oplus 2V \oplus T_1 \oplus T_2$, where $dim_N V=[P:N]-1 $ and $dim_N (T1 \oplus T_2)=[P:N]^2-3[P:N]+1 $ (For more details see \cite{G}). So in this case even without $4$-supertransitivity the above argument holds.
\end{proof}

\subsection{Example: GHJ subfactors for $D_n$ at the trivalent vertex}

We have already seen in \ref{ghj_dn} an example of a noncommuting irreducible quadrilateral whose lower two subfactors have principal graph $A_n$, for each odd value of $n$. The above theorem applies, but the angles can also be computed directly. The computation was done in \cite{GJ} but for convenience we include the argument here. 

\begin{proposition}\label{angle_comp}
Let $N \subset M$ be the GHJ subfactor for $D_n$ at the trivalent vertex, and let $P$ and $Q$ be the pair of intermediate subfactors constructed in \ref{ghj_dn}. Then the angle between $P$ and $Q$ is $cos^{-1} \displaystyle \frac{1}{4cos^2 \frac{\pi}{2n-2}-1}$.
\end{proposition}

\begin{proof}
By $2$-supertransitivity we know that $e_Pe_Qe_P=e_N + \lambda(e_P-e_N)$ for some number $\lambda$ which is the square of the cosine of the angle. 

Let $A_i$, $p$  and $q$ be as in \ref{ghj_dn}. Let $P_0=\mathbb{C}$, $P_1=\{p\}''$, and $P_{i+1}=<P_i,e_i>$ for $i \geq 1$, and similarly for $Q_i$. Then $P=(\cup P_i)'' $, $Q=(\cup Q_i)''$, and the inclusions of towers $P_i \subseteq A_i $, $Q_i \subseteq A_i $ form commuting squares. The element $x=\frac{1}{1-\tau}(p-\tau)$ is in $P$, has norm one, and is orthogonal to $N$. Let $y=\frac{1}{1-\tau}(q-\tau)$, so that $y$ has norm $1$, is in $Q$ and is orthogonal to $N$. Because the squares commute, $E_Q(x)$ must be in $A_1 $ as well as in $Q$ and orthogonal to $N$, and from the Bratteli diagram we see that $E_Q(x)$ must therefore be equal to $cy$ for some constant $c$. We have that $tr(E_Q(x)y)=tr(xy)=tr(\frac{1}{1-\tau}(p-\tau)\frac{1}{1-\tau}(q-\tau)) =\frac{1}{(1-\tau)^2}tr(\tau^2-\tau(p+q))=-(\frac{\tau}{1-\tau})^2$. On the other hand, $tr(E_Q(x)y)=tr(cy \cdot y)=ctr(y^2)=c(\frac{1}{1-\tau})^2tr(q+\tau^2-2\tau q)=c \frac{\tau(1-\tau)}{(1-\tau)^2}=c\frac{\tau}{1-\tau} $. Putting these two equations together we find that $c=-\frac{\tau}{1-\tau} $, and so $E_Q(x)=-\frac{\tau}{1-\tau}y $. Similarly, $E_P(y)=-\frac{\tau}{1-\tau}x$.

Then $E_P E_Q E_P(x)=E_P(E_Q(x))=E_P(-\frac{\tau}{1-\tau}y)=(\frac{\tau}{1-\tau})^2x =\lambda x$, so $\lambda=(\frac{\tau}{1-\tau})^2 $ and the angle is $cos^{-1}\sqrt{\lambda}=cos^{-1}\frac{\tau}{1-\tau}=\cos^{-1}\frac{1}{\tau^{-1}-1} =cos^{-1}\frac{1}{4cos^2 \frac{\pi}{2n-2}-1}$, since the Coxeter number for $D_n$ is $2n-2$.
    
\end{proof}

\section{Subfactors coming from forked Temperley-Lieb algebras}
\subsection{Definition}

We note some properties of the example in the preceeding section.

\begin{lemma}
Let $N \subset M$ be the GHJ subfactor for $D_n$ at the trivalent vertex. Let $p$ and $q$ be as in \ref{ghj_dn}. Then we have: \\
(i) $p$ is orthogonal to $q$.\\
(ii) The sequence of projections $p,e_1,e_2,e_3, ...$ satisfy the Temperley-Lieb relations, with $p$ playing the role of ``$e_0$'' and $\tau=\displaystyle \frac{1}{4cos^2 \frac{\pi}{2n-2} } $. Similary for the sequence $q, e_1, e_2, ...$\\
(iii) $N$ is generated by $e_1,e_2,...$, and $M$ is generated by $p,q,e_1,e_2,...$.
\end{lemma}

\begin{proof}
These are all pretty trivial. (i) follows from the fact that $p$ and $q$ represent distinct vertices at the same level of the Bratteli diagram. (ii) follows from the discussion in \ref{ghj_dn}. (iii) follows since the edges emanting from $p$ and $q$ generate the entire Bratteli diagram.
\end{proof}

This motivates the following definition:

\begin{definition}\label{deffork}
A subfactor $N \subset M$ is called a forked Temperley-Lieb subfactor with parameter $\tau$ if there are projections $p,q,e_1,e_2,e_3,...$ such that:\\
(i) $p$ is orthogonal to $q$.\\
(ii) The sequence of projections $p,e_1,e_2,e_3, ...$ satisfy the Temperley-Lieb relations with parameter $\tau$, and with $p$ playing the role of ``$e_0$''. Similary for the sequence $q, e_1, e_2, ...$\\
(iii) $N$ is generated by $e_1,e_2,...$, and $M$ is generated by $p,q,e_1,e_2,...$.
\end{definition}

Thus the GHJ subfactors for $D_n$ at the trivalent vertex are forked Temperley-Lieb subfactors with parameter $\tau=\displaystyle \frac{1}{4cos^2 \frac{\pi}{2n-2}} $. The forked Temperley-Lieb algebras in the above definition are a special case of algebras associated to $T$-shaped graphs which were studied by Evans and Gould in \cite{EG}. We recall the following result:

\begin{theorem} \label{eg}
\textbf{(Evans-Gould)}
Let $k \geq2$, and let $e_1,e_2,...$ be a sequence of distinct projections satisfying the Temperley-Lieb relations with parameter $\tau > 0$, and let $e_{\bar{k}}$ be another projection. Suppose that we have the following additional relations:\\
(i) $e_ne_{\bar{k}}=e_{\bar{k}}e_n$ if $n \neq k$\\
(ii) $e_{\bar{k}}e_ke_{\bar{k}}=\tau e_{\bar{k}}$ and $e_k e_{\bar{k}} e_k=\tau(1-e_1 \vee ... \vee e_{k-2})e_k$\\
and in the case $k=2$\\
(iii) $e_{\bar{k}}e_1=0$.\\
Let $T_{k,n}$  ($n \in \mathbb{N} \cup \{\infty\}$) be the graph which is constructed by adding one vertex to a Coxeter graph of type $A_n$ and one edge connecting the $k_{th}$ vertex of $A_n$ to the new vertex. Then if $\tau > \frac{1}{||T_{k,\infty}||^2} $ then $\tau=\frac{1}{||T_{k,n}||^2} $ for some $n$. Moreover in this case the $C^*$-algebra generated by $1,e_{\bar{k}},e_1,e_2,...$ is uniquely determined.
\end{theorem}

The case $k=2$ corresponds to the forked Temperley-Lieb algebras (shifting the indices and letting $p=e_1$ and $q=e_{\bar{2}}$). Then $T_{2,\infty}=D_{\infty} $ has norm $2$ and  $T_{k,n}=D_{n+1}$ have norm $2 cos\frac{\pi}{2n} $. Thus if a forked Temperley-Lieb subfactor has parameter $\tau > \frac{1}{4} $, then we must have $\tau=\displaystyle \frac{1}{4cos^2 \frac{\pi}{2n}} $ for some $n$. Also, it follows from the theorem that a forked Temperley-Lieb subfactor is determined by the parameter $\tau $.

\subsection{Classification}

In fact the GHJ subfactors for $D_n$ at the trivalent vertex completely classify noncommuting Type $A_n$ (lower) intermediate subfactors. 

\begin{theorem}\label{bigtheorem}
Let \vpic{quad} {0.6 in}  be a noncommuting, irreducible, hyperfinite quadrilateral with finite index such that $N \subset P$ and $N \subset Q$ are supertransitive with index less than $4$. Then $N \subset M$ is a forked Temperley-Lieb subfactor.
\end{theorem}
\begin{lemma}
Let \vpic{quad} {0.6 in} be a quadrilateral, and let $P \subset M \subset \bar{P}$ be the basic construction. Then $E_{\bar{P}}(e_Q)$ is a scalar multiple of $e_{QP}$.
\end{lemma}
\begin{proof}
There is an ``exchange relation'' for biprojections \cite{Bs3}:\\

\vpic{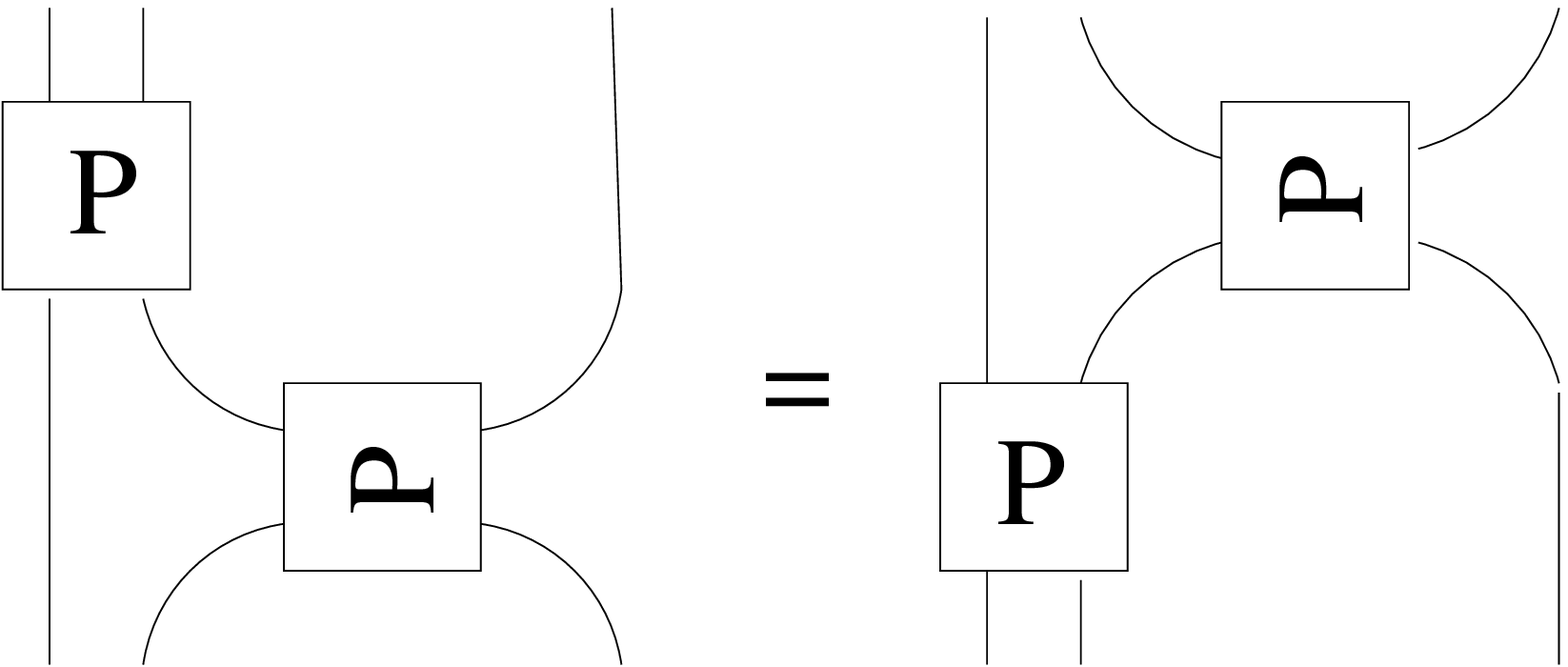} {2.2 in}

which we can apply to obtain:\\

\vpic{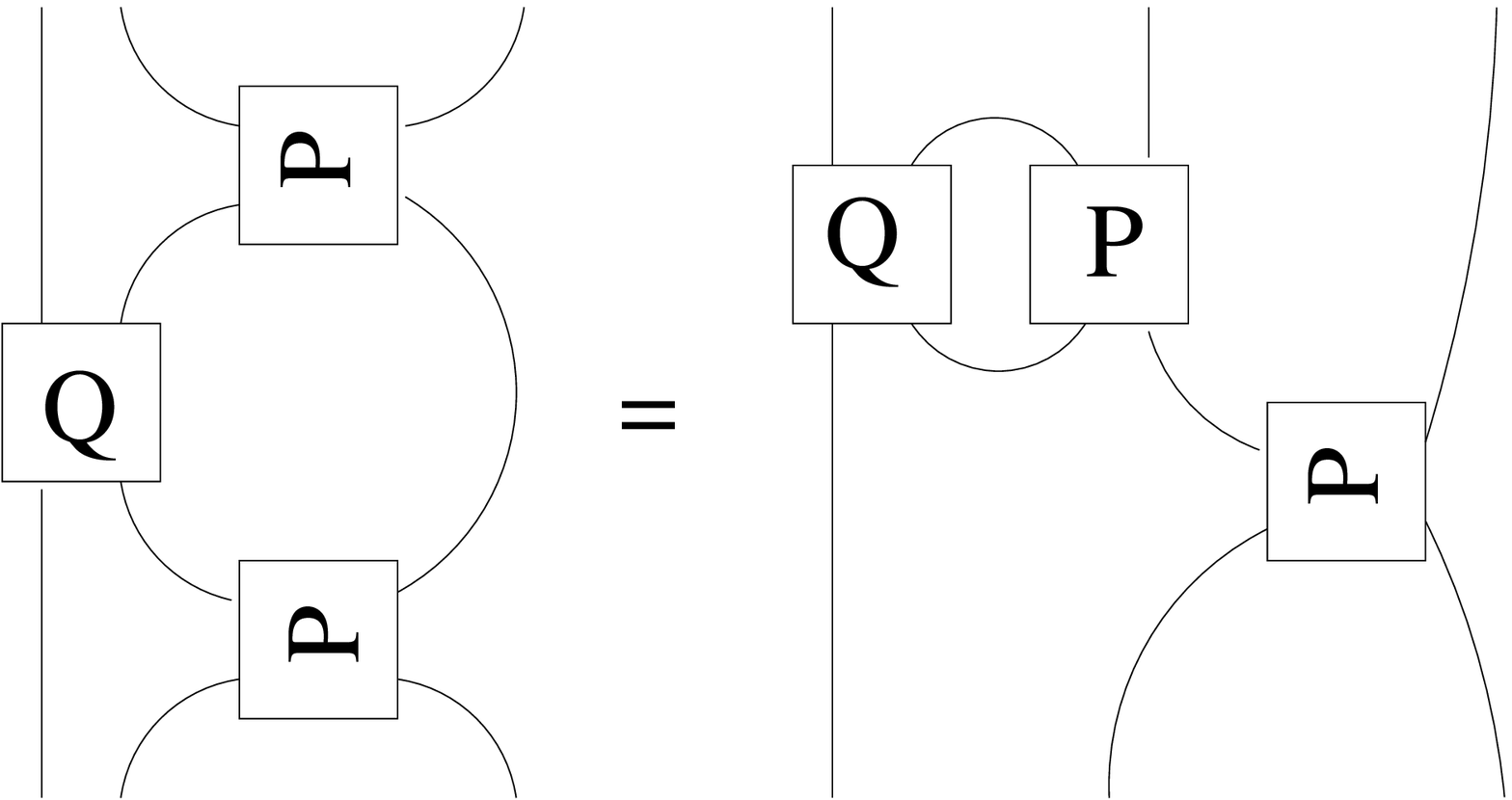} {2.2 in}

Recall that $e_{\bar{P}}$ is a scalar multiple of \vpic{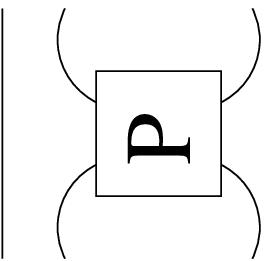} {0.5 in}  . So the left hand side is, up to a scalar, $e_{\bar{P}} e_Q e_{\bar{P}}=E_{\bar{P}}(e_Q)e_{\bar{P}}$, and by \ref{landauproj} the right hand side is, again up to a scalar, $e_{QP} e_{\bar{P}}$.

\end{proof}

\textbf{Proof of \ref{bigtheorem}:}

\begin{proof}
Choose a tunnel $N_{-1} \subset N \subset M$, with corresponding $P_{-1} \subset N \subset P$ and $Q_{-1} \subset N \subset Q$.

Since $N \subset P$-and therefore also $P_{-1} \subset N$- is supertransitive with index less than $4$, it has principal graph $A_n$ for some $n=2,3,4,... $. Then by the uniquess of the hyperfinite $A_n$ subfactor there exist projections $e_1,e_2,...$ satisfying the Temperley-Lieb relations (with $\tau= \displaystyle \frac{1}{4cos^2\frac{\pi}{n+1}}$) such that $N=\{e_1,e_2,e_3,...\}''$ and $ P_{-1}=\{e_2,e_3,e_4,...\}''$. Let $p=e_{P_{-1}} \in P$ be the Jones projection for the inclusion $P_{-1} \subset N$. Then $p$, considered as $e_0$, satisfies the Temperley-Lieb relations along with $e_1,e_2,e_3,...$.

Consider the $N-N$ bimodule map $\phi=\tau^{-1}(e_N+(\tau-1)e_Q):P \rightarrow Q $. By the preceeding lemma $\phi(p)=1-\frac{1-\tau}{\tau}E_Q(p)=1-e_{P_{-1}Q_{-1}}$.

So $q=1-e_{P_{-1}Q_{-1}}$ is a projection orthogonal to $p$, and because $\phi$ is a bimodule intertwiner, it also satisfies the Temperley-Lieb relations (as $e_0$) along with $e_1,e_2,...$. Finally, because $P \vee Q=M$, $p,q,e_1,e_2,...$ generate $M$, and we see that $N \subset M$ is a forked Temperley-Lieb subfactor.
\end{proof}

\begin{corollary}
The series of GHJ subfactors for $D_n$ at the trivalent vertex gives a complete list of noncommuting irreducible quadrilaterals such that $N \subset P$ and $N \subset Q$ are hyperfinite $A_n$ subfactors.
\end{corollary}
\begin{proof}
By \ref{bigtheorem} any such quadrilateral must be a forked TL subfactor, and by \ref{eg} the parameter must be $4cos^{2} \frac{\pi}{2k}$ for some $k$. Since the GHJ subfactors for $D_n$ at the trivalent vertex give examples for each $k$, by ?? it must be isomorphic to one of those GHJ subfactors.  
\end{proof}
\begin{corollary}
The set of numbers which occur as angles between intermediate subfactors subfactors of lower index less than $4$ is exactly \\ $\{ cos^{-1} \displaystyle \frac{1}{4cos^2 \frac{\pi}{2k}-1}|k=3,4,5,... \}$.
\end{corollary}
\begin{proof}
By \ref{angle_thm} if $[P:N],[Q:N]<4 $ then the angle is equal to $\cos^{-1} \displaystyle \frac{1}{[P:N]-1}$. On the other hand, we know from \cite{J3} that the only admissible index values less than $4$ are $4cos^2 \frac{\pi}{k}$, $k=3,4,5,... $. So the only possible angle values for quadrilaterals with small lower indices are given by the series $cos^{-1} \frac{1}{4cos^2 \frac{\pi}{k}-1} $. By \ref{angle_comp} the GHJ subfactors for $D_n$ at the trivalent vertex give examples of all these values for even $k $. To achieve the values for odd $k$, $N \subset P $ and $N \subset Q $ would have to have as principal graphs Dynkin diagrams with odd Coxeter numbers, of which the  only candidates are $A_{2k}$. But by \ref{bigtheorem} this would mean that the quadrilateral is a forked TL subfactor, and by \ref{eg} this is impossible for an odd Coxeter number.  
\end{proof}

%


\bibliographystyle{amsplain}

\5
\5
\end{document}